\documentclass[11pt]{amsart}
\usepackage{a4wide}
\usepackage{amsthm,amsfonts,amssymb}
\usepackage{amsrefs}
 \usepackage{graphicx}
\usepackage{hyperref}

\begin{document}
\newcommand{\s}{\sigma}
\renewcommand{\l}{\lambda}
\newcommand{\F}{\Phi}
\renewcommand{\t}{\theta}
\newcommand{\C}{\mathcal C}

\def\RR{\mathbb R}
\def\ZZ{\mathbb Z}
\def\NN{\mathbb N}
\def \cf {{\mathcal F }}
\def \Q {{\mathcal Q }}
\def \a {\alpha }
\def \b {\beta }
\def \f {\varphi }
\def \L {\Lambda } 
\def \o {\omega }
\def \O {\Omega } 
\def \D {\Delta }
\def \d {\delta }
\def \g {\gamma }
\def \e {\epsilon }
\def \r {\rho }

\newtheorem{T}{Theorem}
\newtheorem*{theorem}{Theorem}
\newtheorem{lemma}{Lemma}
\newtheorem{Corollary}{Corollary}
\newtheorem{maintheorem}{Theorem}
\newtheorem{Proposition}{Proposition}
\newtheorem{Remark}{Remark}
\newtheorem{Definition}{Definition}
\theoremstyle{definition}
\newtheorem{remark}{Remark}

\title[Hyperbolicity 
of horseshoes with internal tangencies]
 {Some non-hyperbolic systems with 
strictly non-zero Lyapunov exponents for all invariant measures: 
\\ Horseshoes with internal
 tangencies}

\author{Yongluo Cao}
\address{Department of Mathematics, Suzhou University, Suzhou 215006, 
Jiangsu, P.R. China}
\email{ylcao@suda.edu.cn}
\author{Stefano Luzzatto}
\address{Dept. of Mathematics, Imperial College London, 
UK}
\email{stefano.luzzatto@ic.ac.uk}
\urladdr{http://www.ma.ic.ac.uk/\textasciitilde 
luzzatto}

\author{Isabel Rios}
\address{Instituto de 
Matem\'atica,
Universidade Federal Fluminense 
(UFF),
Brazil}
\email{rios@mat.uff.br}
\thanks{This work was partially supported by CAPES and FAPERJ
(Brazil). The authors acknowledge the hospitality of Imperial College 
London where this work was carried out.}

\begin{abstract} 
We study the hyperbolicity of a
class of horseshoes 
exhibiting an \emph{internal} tangency, i.e.
a point of homoclinic 
tangency accumulated by periodic points. In particular these systems 
are strictly \emph{not} uniformly hyperbolic. However we show that 
 all the Lyapunov exponents of all invariant measures are uniformly 
 bounded away from 0. This is the first known example of this kind. 
\end{abstract}  
\date{9 September 2004}


\maketitle

\section{Introduction}\label{int}

\subsection{Hyperbolicity and tangencies}
\begin{figure}[htb]
\begin{center}
    \includegraphics[width=6cm]{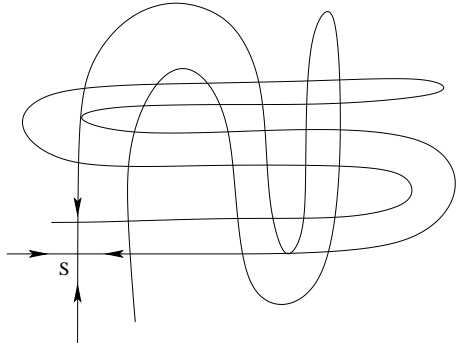}
 \end{center}   
\caption{\label{map} 
Homoclinic tangencies inside the limit set}
\end{figure}

We consider $C^2$ diffeomorphisms 
on Riemannian surfaces.  
Our goal is to study the hyperbolic 
properties of a class of  maps exhibiting a homoclinic tangency
associated 
to a fixed saddle point $S$, as in Figure 
\ref{map}.
We assume without loss of generality that we are 
working on \( 
\mathbb R^{2} \) and in the standard Euclidean norm.
We recall that a compact invariant set 
\( 
\Lambda \) is \textbf{uniformly hyperbolic}
    if there exist constants 
   \( C>0,  \quad  \sigma > 1 > \lambda > 0 \)  and a continuous, \( 
D\Phi \)-invariant,  
   decomposition 
   \(T_{x}\Lambda = 
E^{s}_{x}\oplus E^{u}_{x} \) of the 
   tangent bundle over \( \Lambda \) such 
that for all \( x \in 
\Lambda \) and all \( n\geq 1 \) we have 
 \begin{equation} \label{UH}
\|D\Phi^{n}|_{E^{s}_{x}}\| \leq C\lambda^{n} 
\quad\text{ and } \quad
\|D\Phi^{n}|_{E^{u}_{x}}\| \geq C^{-1}\sigma^{n}. 
\end{equation}
By standard hyperbolic theory, every point \( x \) in \( \Lambda \) 
has stable and unstable manifolds \( W^{s}_{x}, W^{u}_{x} \) 
tangent to the subspaces \( E^{s}_{x} \) and 
\( E^{u}_{x} \) respectively, and thus in particular 
\emph{transversal} to each other. 
The presence of the tangency therefore implies that the 
dynamics on \(\Lambda \) cannot be uniformly hyperbolic.  

We emphasize at this point that, in the case we are considering, 
the homoclinic tangency is accumulated by 
transverse homoclinic orbits which in turn are accumulated by 
periodic points. Thus it constitutes an intrinsic obstruction to 
uniform hyperbolicity which cannot be resolved 
by simply ignoring the orbit of tangency. Most of the 
classical theory of homoclinic bifurcations for diffeomorphisms 
(see \cite{PalTak93} and references therein) 
considers the unfolding of homoclinic tangencies 
 \emph{external} to the set \( \Lambda \) to which they are 
 associated, thus causing no real issues with the hyperbolicity at the 
bifurcation parameter. The main goal of such a theory has often 
been to study the hyperbolicity and the occurrence of tangencies 
in a neighbourhood of the orbit 
of tangency \emph{after} the bifurcation. The presence of an internal 
tangency gives rise to a more subtle situation and it has only recently been 
shown that this can actually occur as a \emph{first bifurcation} 
\cites{Kir96a, Asa97, CaoKir00, BedSmi, Rio01}. Part of the motivation 
of the present paper is to study the global dynamics and hyperbolicity
\emph{at this bifurcation parameter}. 

\subsection{Hyperbolicity and Lyapunov exponents}
Another part of the motivation for this result is to give an example 
of a compact invariant set which is as ``uniformly'' hyperbolic as 
possible in the ergodic theory sense, but still not uniformly 
hyperbolic. To formulate this result precisely,  
let \( \mathcal M (\Phi) \) denote the set of all 
\( \Phi \)-invariant probability measures \( \mu \) on \( \Lambda \). By 
the classical Multiplicative Ergodic Theorem of Oseledet's there is a 
well defined set \( \mathcal L(\mu) \) of \emph{Lyapunov exponents} 
associated to the measure \( \mu \), we give the precise definitions 
below. We say that 
    the measure \( \mu\in \mathcal M(\Phi) \) is \emph{hyperbolic} if 
    all the Lyapunov exponents are non-zero. 
The existence of an invariant measure with non-zero Lyapunov exponents 
indicates a minimum degree of hyperbolicity in the system. 
A stronger requirement  is that all invariant measures 
    \( \mu \) are hyperbolic 
    and of course an even stronger requirement is that they are all 
    ``uniformly'' hyperbolic in the sense that all Lyapunov exponents 
    are uniformly bounded away from 0. 
This condition is  clearly satisfied  for uniformly hyperbolic systems but, as 
we show in this paper, it is strictly weaker. 

The class of examples we are 
interested in were first introduced in
\cite{Rio01} and constitute perhaps the simplest 
situation in which an internal tangency can occur as a first 
bifurcation. In section \ref{maps} we give 
the precise definition of this class.
For this class we shall then prove the following
\begin{theorem}
All Lyapunov exponents of all measures in  \( \mathcal 
 M(\Phi) \) are uniformly bounded away from zero.
\end{theorem}

As an immediate corollary we have the following statement which is in 
itself non-trivial and already remarkable. 

\begin{Corollary}
\( \Phi\) is uniformly
hyperbolic on periodic 
points.
\end{Corollary}

We recall that uniform hyperbolicity on periodic points means that
there exists constants 
\( \sigma > 1 > \lambda > 0 \)  
   such 
that  for each 
   periodic point \( p \in \Lambda \) of period \( 
k \), the 
   derivative \( D\Phi^{k}_{p} \) has two distinct real 
eigenvalues 
   \( \tilde \sigma, \tilde\lambda \) with 
\(
   |\tilde 
\sigma| > \sigma^{k} > 1 > \lambda^{k} >     |\tilde\lambda| 
> 0. 
\)
   Notice that the bounds for the eigenvalues are exponential in \( k \).

As far as we know this is the first known example of this kind, 
although it is possible, and indeed even likely, that such a property 
should hold for more complex examples such as Benedicks-Carleson and 
Mora-Viana parameters in H\'enon-like families \cites{BenCar91, 
MorVia93} and horseshoes at the 
boundary of the uniform hyperbolicity domain of the H\'enon family. 
The weaker result on the uniform hyperbolicity of periodic points 
has been proved recently for both cases in \cite{WanYou01} and 
 \cite{BedSmi} respectively. .  
Other known examples of non-uniformly hyperbolic diffeomorphisms include
cases in which the lack of uniformity comes from the presence of
``neutral'' fixed or periodic points. In these cases, the Dirac-\( 
\delta \) measures on such periodic orbits are invariant and have a 
zero Lyapunov exponent. 

It is interesting to view our result in the light of some recent work 
which appears to go in the opposite direction: if a 
compact invariant set $\Lambda$ admits an invariant splitting  \(T_{x}\Lambda = 
E^{1}_{x}\oplus E^{2}_{x} \) such that, for a total measure set of points $x\in \Lambda$, 
the Lyapunov exponents are positive in $E^{1}_{x}$ and negative in $E^{2}_{x}$, then 
$\Lambda$ is uniformly hyperbolic, \cites{AlvAraSau03, Cao03, CaoLuzRioHyp}.
Here, the Lyapunov exponents are not even required to be uniformly bounded away 
from zero. Thus the existence of at least one orbit as which the 
splitting degenerates is a necessary condition for a situation such 
as the one we are considering, in which the Lypaunov exponents are all 
non-zero but \( \Lambda \) is strictly not uniformly hyperbolic.

The concept of uniform hyperbolicity of the periodic 
points and of measures plays an important role
in the general theory of one-dimensional dynamics. 
In some situations, such as for certainclasses of 
smooth non-uniformly expanding 
unimodal maps these notions have been been shown to be 
equivalent to each other and to various other properties usually associated to uniform 
hyperbolicity such as exponential decay of correlations 
\cite{NowSan98}.
Higher dimensional cases are generally much harder  due to the 
complexity of the geometrical and dynamical structure of the examples 
and progress is only starting to be made.

\subsection{Overview of the paper}
After defining the class of systems of interest in section 
\ref{maps}, 
in section \ref{splitting} we construct a field of cones which are 
invariant under the derivative map. In Section 
\ref{periodicpoints}
we show that some \textbf{non-uniform hyperbolicity}  is 
satisfied off the orbit of tangency: the expanding and contracting
directions of points which are very 
close to the orbit of tangency are
almost aligned and therefore the 
contraction and/or 
expansion may take arbitrarily long time to become 
effective. We formalize this by saying that there exists a constant \( 
C_{x} >0\), \emph{depending on the point} \( x \) (compare with \eqref{UH}),  
such that for all \( x \) not in the orbit of tangency we have 
 \begin{equation} \label{NUH}
\|D\Phi^{n}|_{E^{s}_{x}}\| \leq C_{x}\lambda^{n} 
\quad\text{ and } \quad
\|D\Phi^{n}|_{E^{u}_{x}}\| \geq C^{-1}_{x}\sigma^{n}. 
\end{equation}
with \( 
C_{x} \) arbitrarily small near the points of the orbit of 
tangency, 
i.e. 
\(\inf_{x\in\Lambda}C_{x}=0. \)
One way to show that the 
non-uniformity is not too extreme is to show that there is 
 a  \emph{uniformly hyperbolic core} in the 
following sense: there exists a region $W$ containing one point of the orbit 
of tangency, and there exists constants 
\( \sigma > 1 > \lambda > 0 \)  such 
that,  for each point $p\in\Lambda \cap W$ that returns to $W$ at the $k$-th iterate, we have
\begin{equation} \label{UHk}
\|D\Phi^{k}|_{E^{s}_{p}}\| \leq \lambda^{n} 
\quad\text{ and } \quad
\|D\Phi^{k}|_{E^{u}_{p}}\| \geq \sigma^{n}. 
\end{equation}
In fact, we show that this uniform hyperbolicity occurs well before the return of the 
orbit of $p$ to $W$. Those exponential estimates are valid from the moment that 
the orbit of $p$ leaves a certain neighborhood of the fixed point $S$, and they keep 
working as long as the orbit wanders around $\Lambda \setminus W$.
Finally, in section \ref{lyapunov} we show that this implies the 
statement of our main theorem.

\section{Horseshoes with internal 
tangencies}\label{maps}

We start with a geometrical definition of the class of maps 
under 
consideration. 
  Let $\F:\mathbb R^{2}\to \mathbb R^{2}$ be a  
\( C^{2} \) 
  diffeomorphism and 
let \( Q = [0,1]\times [0,1] \) 
denote the unit square in \( \mathbb R^{2} \). 
We suppose that there exist 
5 ``horizontal'' regions in \( Q \) which 
are mapped as in Figure 
\ref{the map}: $R_i^,=\F(R_i)$ for $i=1,\ldots ,5$.
\begin{figure}[htb]
\centerline{\includegraphics[width=9cm]{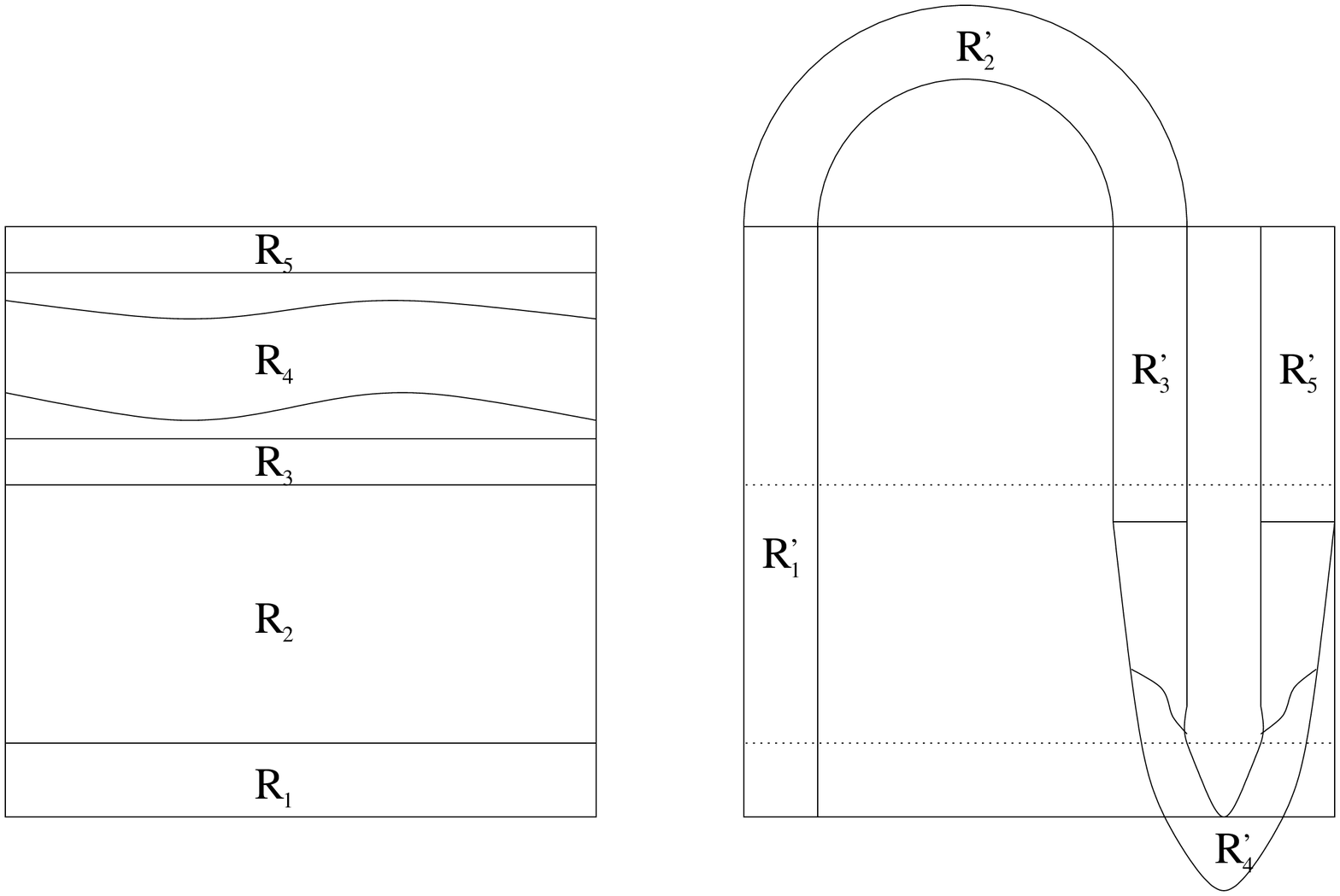}}
\caption{\label{the 
map}The map $\F_c$}
\end{figure}
We suppose that regions \( R_{1}, R_{3}, 
R_{5} \) are mapped affinely 
to their images, with derivative 
$$D\F(x,y)=\left(\begin{array}{cc}
\l & 0 \\ 0 & \sigma \end{array} \right) 
\quad\text{and}\quad 
D\F(x,y)=\left(\begin{array}{cc}
- \l & 0 \\ 0 & 
-\sigma \end{array} \right) 
$$
in \( R_{1}\cup R_{3} \) and \( R_{5} \) 
respectively, for two 
constants \( \lambda, \sigma \) satisfying
\( 
\sigma > 3 \)  and \( \lambda < 1/3 \). More explicitly 
we suppose that 
$\F(x,y)=(\l x,\sigma y)$ in 
in \( R_{1} \) ($0\leq y\leq 
\sigma^{-1}$), 
\( \F(x,y)=(\l x+(1-\l),\sigma y-(\sigma-1)) \) in \( R_{5} \) 
($1- 2/3\sigma\leq
 y\leq 1$) and that \( R_{3}' \) is a vertical strip 
parallel to \( R_{5}' \). 
 In particular there is a hyperbolic fixed 
point at the origin \( (0,0) 
 \) whose stable and unstable manifolds 
contain the lower and left 
 hand side of the square \( Q \) respectively. 

We emphasize that the explicit nature of the map in these regions is 
for simplicity only, and to allow us to concentrate on the strategy 
for dealing with the tangency. This could be weakened significantly, 
for 
example by assuming only some uniformly hyperbolic structure in 
regions \( R_{1}, R_{3}, R_{5} \).  
 Region \( R_{2} \) is mapped outside 
the square \( Q \) and thus we 
 do not need to make any particular 
assumption on its form.

 Region 4 is mapped to a ``fold'' \( R'_{4} \) 
which contains a point 
of the orbit of tangency. Here we make 
some
``non-degeneracy assumptions'' essentially stating that vertical lines 
in \( 
R_{4} \) are mapped to non-degenerate parabolas. 
 More precisely, we 
suppose that there is a region \( \hat R_{4}\subset 
  R_{4} \) bounded 
by two disjoint curves
 so that points in $R_4\setminus \hat R_{4}$, 
are
mapped inside region $R_2$ with second coordinate 
greater than 
$\sigma^{-1}$. Then, for each $x_0$ in $[0,1]$ 
we have that  $\F 
(\{(x_0,y):y\in 
{\mathbb R}\}\cap \hat R_4) $
is contained in the graph of the 
map  
\begin{equation}\label{parabolas}
f_{x_0}(x)=c(x-q)^2-\l x_0
\end{equation}
for some fixed \( c>0 \) and some 
$q\in (2/3,1)$.  
Thus all vertical lines in \( \hat R_{4} \) are 
mapped to parabolas 
with constant curvature \( c \). Again, this could be 
weakened 
significantly although we do assume for less trivial reasons 
that the 
curvature \( c \) is sufficiently large in relation to the 
constants 
\( \lambda \) and \( \sigma \). 
Notice that the point \( (q,0) 
\) is a point of tangency of the 
stable and unstable manifolds of the 
fixed point at the origin. 
To control the global properties of the 
family of parabolas we also 
assume that 
 \[
 \left\|{\partial \F \over 
\partial y} (0,y)\right\|\geq \sigma 
 \text{ and } 
 \left\|{\partial 
\F \over \partial x}(\F^{-1}(q,0))\right\|= \l
\]
and that, for every \( 
(x,y)\in\hat R_{4} \), 
$$\left< {\partial \F \over \partial y} 
(x,y),
{\partial \F \over \partial x} (x,y) \right>=0$$
A completely explicit 
example of a map 
\( \Phi \) satisfying all these conditions is 
given 
in \cite{Rio01}. We now define the limit set 
\[ 
\Lambda = 
\bigcap_{n\in\mathbb Z} \Phi^{n}(Q).
\]
It follows from the construction that \( 
\Lambda \) is non empty and 
contains (at least) one orbit of tangency 
between stable and unstable 
manifolds. Our assumptions make it 
easy to control the dynamics 
in regions \( R_{1}, R_{3}, R_{5} \) and 
therefore our main objective 
is to control the dynamics of points 
returning to the fold \( \hat R'_{4} \) 
and to a neighbourhood of the 
tangency. 
Thus we concentrate on the 
\emph{first return map} from \(  \hat 
R'_{4} \) to itself and show 
that this map satisfies some strong 
hyperbolicity estimates from 
which we can deduce our main estimate on the 
hyperbolicity of 
periodic points.

\section{The  conefield}\label{splitting}

Consider the foliation ${\cf}$, of $\F(Q)$, whose leaves are 
images of
 vertical lines by the map $\F$ (in $\hat R'_4$, the leaves 
are parabolas). To each point $(x,y)$ we
associate the tangent 
direction, $l(x,y)$, to the leaf of ${\cf}$
which passes through $(x,y)$, at 
this point ($l(x,y)$ 
is parallel to ${\partial \F \over \partial 
x}(\F^{-1}(x,y))$.
Let $a(x,y)$ be the angle between $l(x,y)$ and the 
horizontal. 
There is only one point for which \( a=0 \) and this is the point of 
tangency \( (q,0) \). At all other points we have 
with $0< a(x,y)\leq 
\pi /2$. Then, we divide \( \Lambda \) into 3 
regions, according to the 
value of \( a \): 
$$L=\{a =\pi/2\}, \ R=\{ \pi/3< a <\pi/2\}, 
V=\{ 0 
< a \leq\pi/3\}.$$
Notice that 
\[ 
\L\cap (R'_1\cup R'_3\cup R'_5) 
= L \text{ and } 
V\cup R = \L\cap(R'_4\setminus (q,0) ) 
\]
 
For each point $P \in R\cup V$ we 
define a \emph{cone} \( \mathcal C(P) \) 
 of vectors in 
the tangent 
space \( T_{P}\mathbb R^{2} \) as follows: 
let $P=(\xi ,\eta)$, and $P_{i}=\F^i(P)=(\xi_i, \eta_i)$ for any integer $i$. 
All cones will be centered on the vertical and, for convenience, we 
define a \emph{standard cone} 
\begin{equation}\label{standard}
{\C}=\left\{(v_1,v_2)\in {\mathbb R}^2: 
{{|v_1|}\over{|v_2|}} \leq 
\sqrt{3} \right\}. 
\end{equation}
For $P$ in $ R$ we simply let \( \mathcal C(P) = \mathcal C\). 
For $P\in V$, we define cones which extend to each side of the 
vertical by three times the (cotangent of the) angle of the line \( l(P) 
\) defined above. More precisely, recall that \( l(P) \) is tangent to 
the graph of \( 
f(x)=c(x-q)^{2}-\xi_{-1}\lambda \), see \eqref{parabolas}, and 
therefore \( l(P)=\{ (v_1,v_2)\in \mathbb R^2: 
v_2=\pm 2 v_1\sqrt{c(\eta+\l \xi_{-1})}
\). Then we let
\[
\mathcal C (P)=
\left\{(v_1,v_2)\in \mathbb R^2:
\frac{|v_1|}{|v_2|}\leq 
\frac{3}
{2\sqrt{c(\eta+\l \xi_{-1})}}\right\}
\]
Before defining the cones in the remaining points 
in \( L \) we show that the conefield defined above in \( R\cup V \) 
is invariant under the first return map to \( R\cup V \). 
As a first step towards proving this invariance we estimate the first 
return time. 

\begin{lemma}
\label{return}
Consider 
$P=(\xi,\eta) \in R\cup V$, and $n$ the smallest
positive integer 
such that $P_n 
\in R\cup V$. Then  
\begin{equation}
\label{estim1n}
n\geq\max\left\{{{\ln{1/\eta}}\over{\ln{\sigma}}},
{{\ln{1/\xi_{n-1}}}\over{\ln{\l^{-1}}}}\right\}.
\end{equation}
\end{lemma}

\begin{proof}
Note 
that, since $V\subset R_1$, we have $\F^i(P)=
(\l^i\xi ,\sigma^i 
\eta)$, provided that $\sigma^{i-1}\eta \leq
\sigma^{-1}$. 
In fact, the images $P_i$ of the point 
$P$ will stay in region $R_1$ 
while
$i$ is such that 
$\eta\sigma^{-i}<\sigma^{-1}$. Since $P_{n-1}\in R_4$, there must be at 
least one more 
linear iterate of $P_i$ before its orbit visits the pre-image of $R\cup 
V$. 
So, we have $\eta\sigma^{n-1}> \sigma^{-1}$. 
Analogously, 
using the inverse map, we have 
\(
n\geq 
 {{\ln{1/\xi_{n-1}}}\over{\ln{\l^{-1}}}}
\).
\end{proof}

We shall show below that these estimates mean that cones are 
sufficiently contracted before returning,  guaranteeing that they are 
mapped back strictly into existing cones even if they started out 
very wide. In particular we shall use the following simple
\begin{Corollary}
\label{returncor}
Consider 
$P=(\xi,\eta) \in R\cup V$, and $n$ the smaller
positive integer 
such that $P_n 
\in R\cup V$. Then 
for
this $n$ we 
have
\begin{equation}
\label{lambdan1}
{{\l^n}\over{\sigma^n}}<
\eta^{1-{{\ln{\l}}\over{\ln{\sigma}}}}
\end{equation}
and
\begin{equation}
\label{lambdan22}
{{\l^n}\over{\sigma^n}}<(\l 
\xi_{n-1})^{1-{{\ln{\sigma}}\over{\ln{\l}}}}\,.
\end{equation}
\end{Corollary}

\begin{proof}
The first inequality follows from \eqref{estim1n} which gives
    \begin{equation*}
{{\l^n}\over{\sigma^n}}<\left({{\l}\over{\sigma}}\right)^{{\ln{1/\eta}}
\over{\ln{\sigma}}}=\left({{1}\over{\eta}}\right)^{{\ln{\l}}\over{\ln{\sigma}}}
\left({{1}\over{\sigma}}\right)^{{\ln{1/\eta}}\over{\ln{\sigma}}}=
\eta^{1-{{\ln{\l}}\over{\ln{\sigma}}}}.
\end{equation*}
The second follows also by \eqref{estim1n} by a similar straightforward 
calculation. 

\end{proof}
We are now ready to prove the invariance of the conefield for the 
first return map. 
\begin{lemma}
\label{inclusion}
There exists $c_1>0$ such 
that, if  $c>c_1$ in the definition of $\F$, 
$P\in V\cup 
R$ and $P_n$ 
is the first
positive iterate of $P$ in $R\cup V$, 
then
$$D\F^n(P)({\C}(P))\subset {\C}(P_n)\,.$$
\end{lemma}

\begin{proof}
Since $\C(P)$ 
is
centered in the vertical line, we have, 
after
applying
$$D\F^{n-1}(P)=
\left(\begin{array}{cc} \pm \l^{n-1} & 0 \\ 
0 & \pm
\sigma^{n-1}\end{array}\right)$$
to the vectors of ${\C}(P)$, 
that
$D\F^n(P)\left({\C}(P)\right)$ is a cone centered in the vertical 
line 
at
the point $P_{n-1}$ whose
width is 
${\l^{n-1}}\over{\sigma^{n-1}}$ times the width of ${\C}(P)$. 
After
applying $D\F$ to this
cone, 
since $P_{n-1}\in R_4$, we obtain a cone centered 
in
$l\left(P_{n}\right)$ with width
smaller than ${{\l^n}\over{\sigma^n}}$ times the original 
width.

There are many cases to be considered 
depending on the location of $P$ and 
$P_n$ in $R\cup V$.
Here we present the 
case where $P$ and $P_n$ are both in $V$. The 
other cases are made 
following the same steps. 
Let $P$ and $P_n$ be in 
$V$. Then 
$${\C}(P)=\left\{(v_1,v_2)\in{\mathbb R}^2:{ 
{|v_1|}\over{|v_2|}}
\leq {{3}\over{2\sqrt{c(\eta +\l 
\xi_{-1})}}}\right\}$$
and
$${\C}(P_n)=\left\{(v_1,v_2)\in{\mathbb R}^2:{{|v_1|}\over{|v_2|}}
\leq 
{{3}\over{2\sqrt{c(\eta _n+\l \xi_{n-1})}}}\right\}\,.$$
Since $0<a(P_n)<\pi /3$, it 
satisfies
$$\arctan{\left({{\tan{a(P_n)}}\over {3}}\right)}< {{a(P_n)}\over 
{2}}\,.$$
This implies that 
a cone centered in $l(P_n)$ with width ${{2\sqrt
{c(\eta _n+\l \xi_{n-1})}}/ 
{3}}$ will be properly contained in  ${\C}(P_n)$.
Then, in order to 
obtain
\(D\F^n(P)({\C}(P))\subset {\C}(P_n) \)
we need
\begin{equation*}
{{\l^n}\over{\sigma^n}}{{3}\over{2\sqrt{c(\eta +\l 
\xi_{-1})}}}<
{{2\sqrt{c(\eta _n+\l \xi_{n-1})}}\over {3}}
\end{equation*}
or, equivalently, 
\begin{equation}
\label{lambdan3}
{{\l^n}\over{\sigma^n}}<
{{4c\sqrt{(\eta_n+\l \xi_{n-1})}\sqrt{(\eta +\l \xi_{-1})}}\over 
{9}}\,.
\end{equation}
We now distinguish two cases. 
In the case where $\eta\leq \l \xi_{n-1}$ 
we have 
\(
\eta< \sqrt{(\eta_n+\l \xi_{n-1})}\sqrt{(\eta +\l \xi_{-1})}
\)
and therefore it is enough to show that \( \lambda^{n}/\sigma^{n}< 
4c\eta/9 \).
By \eqref{lambdan1} this follows as long as we 
 have 
\begin{equation}
\label{novo}
\eta^{1-{{\ln{\l}}\over{\ln{\sigma}}}}< 
{{4c\eta}\over{9}}\,,
\end{equation}
for all $\eta \leq 1$ (remember that we are working with 
points of $[0,1]\times[0,1]$). Since
$\ln{\l}/\ln{\sigma}$ is 
fixed
and negative, the
condition \eqref{novo} holds if $c$ 
is big enough. 

In the case where $\eta> \l \xi_{n-1}$ we argue in the same way and, 
using \eqref{lambdan22}, 
reduce the problem to showing that
\begin{equation*}
(\l \xi_{n-1})^{1-{{\ln{\l}}\over{\ln{\sigma}}}}< 
{{4c\l\xi_{n-1}}\over{9}}\,,
\end{equation*}
for any $\xi_{n-1}\leq 1$. Again this follows as long as \( c \) is 
sufficiently large. 
The other cases follow 
analogously. 
\end{proof}

We can now 
extend the conefield to the set $\L\setminus \{ 
\F^i(q,0), i
\in 
\ZZ \}$, in a natural way 
by considering the images of of all cones in \( R\cup V \) and taking 
slightly wider cones at each point. 
In this way we obtain a conefield defined at every point outside the 
orbit of tangency such that
$$D\F(P)({\C}(P))\subset 
{\C}(\F(P))\,.$$
for every point \( P \). Notice that these cones can be arbitrarily 
wide close to the orbit of tangency. 
For points \( P \) which do not enter \( R\cup V \)  in either forward or 
backward time, we simply let \( \mathcal C(P) \) be the standard cone 
\( \mathcal C \), see \eqref{standard}. For points \( P \) which intersect \( 
R\cup V \) only in forward time, we define 
\[ 
C(P) = \mathcal C \cap D\Phi^{-i}(\mathcal C(\Phi(P)))
\]
where \( i >  0 \) is the first time for which \( \Phi(P)\in R\cup V \) 
and \( \mathcal C\) is that standard cone. 

The stable 
cone
field is defined assigning to each
$P\in 
\L\setminus \{ \F^i(q,0), i \in \ZZ \}$, the closure of the 
complement
of 
${\C} (P)$, and satisfies the inclusion condition for the inverse 
map.

\section{Uniform hyperbolicity for the escape time}
\label{periodicpoints}

We fix a neighbourhood $W$ of the tangency point 
$(q,0)$ of radius $1/c$. Let $P=(\xi, \eta)\in W$. 
We start with a simple estimate of the number of iterations 
$n$ it takes before  $P_n$ falls 
outside the set  $[0,1]\times[0,1/3]$ ($n$ is the escape time of $P$). 

\begin{lemma} 
\label{pontos}
We have \( n\geq -\log 3\eta/\log\sigma \). 
\end{lemma}

\begin{proof} 
Note that
$W$ is 
contained in the
domain where $\Phi$ is linear, and, as long as the 
second coordinate of 
the
image of a point of $\Lambda$ is less 
than
$1/3$, it is contained in $R_1$, the second coordinate of $P$ 
being
multiplied by $\sigma$ at each
iteration. 
\end{proof}

For $P=(\xi,\eta)$ as above, let $v$ be a vector 
contained in ${\C}(P)$ and $v^n=D\Phi^n_P v$. 
\begin{lemma} 
\label{vetores}
We have
$\Vert v^n \Vert \geq 
\sigma^n|v_2|.$
\end{lemma}
\begin{proof}
As we saw in the last section, since $\Phi^n$ is 
linear in $W$, we have
$$\Vert v^n \Vert =\Vert (\lambda^n v_1, 
\sigma^n v_2) \Vert \geq
\sigma^n|v_2|.$$

\end{proof}

\begin{lemma} 
\label{expansao}With the notation above,  if $c$ is sufficiently big, then
$$\frac{\Vert v^n 
\Vert}{\Vert v \Vert}\geq \sigma^{n/2}$$
\end{lemma}

\begin{proof}
Since $v 
\in {\C}(P)$, we have that 

\begin{equation}
\label{tamanho}
\Vert v 
\Vert = \sqrt{v_1^2+v_2^2}\leq 
\sqrt{\left(\frac{3}{2\sqrt{c\eta}}\right)^2+1}|v_2|
\end{equation}
Using lemmas \ref{vetores} and \ref{pontos}, we 
obtain
$$\frac{\Vert v^n \Vert}{\Vert v 
\Vert}\geq
\frac{\sigma^n}{\sqrt{\frac{9}{4c\eta}+1}}=
\frac{\sigma^n2\sqrt{c\eta}}{\sqrt{9+4c\eta}}.$$
By lemma 
\ref{pontos} we have that $\sigma^{n/2}\geq \frac{1}{\sqrt{3\eta}}$ 
and therefore 
$$\frac{\Vert v^n \Vert}{\Vert v \Vert}\geq
\frac{2 
\sqrt{c}\sigma^{n/2}}{\sqrt{3}\sqrt{9+4c\eta}}.$$
Since 
$|\eta|<1/c$, we can 
write
$$\frac{\Vert v^n \Vert}{\Vert v \Vert}\geq 
\frac{2\sqrt{c}\sigma^{n/2}}{\sqrt{39}}.$$
The result follows if $c$ is sufficiently large (greater than 
$\sqrt{39/4}$ in this case).
\end{proof}

\section{Lyapunov exponents}
\label{lyapunov}

Let $W_j$ be the connected 
component of $\F^j(W )\cap R_1$ containing 
$\F^j(q,0)$, and 
$\tilde{W}=\cup_{j\geq 0}W_j$. In this section we estimate the growth of the unstable vectors of points in
$\L$ outside $\tilde{W}$, under the action of $D\F$. We also compute bounds for the Lyapunov exponents of $\F$, and prove the main result.

\begin{lemma}\label{first}
There exists $\sigma_1>1$ such that, if $P \in\L \setminus \tilde{W}$, and $v\in \C(P)$, then $\Vert D\F_Pv\Vert > \sigma_1\Vert v \Vert$. 
\end{lemma}
\begin{proof}
First we claim that, if $P \in\L \setminus \tilde{W}$,  
we have \( \C(P)\subset \mathcal C \). 
Indeed, by construction, if $P\in R_4^,$ and
at distance at least $1/c$ from the tangency point 
$(q,0)$, then we 
have that 
$\eta +\l \xi_{-1} \geq 1/c$, and the width of 
${\C}(P)$ is less 
than $3/2 < \sqrt{3}$. 
Notice that vectors in  \( \C \) 
grow by at least a constant 
factor $\sigma_1>1$ at each 
iteration where $\sigma_1>1$ is the rate of 
growth of the vector 
$(\sqrt{3},1)$ by the linear map $(x,y)\mapsto (\l x, \sigma y)$).
If $P$ is a point outside the set $\tilde{W}$, then we have already 
\( \C(P)= \mathcal C \), and the estimate applies to the unstable vectors of $P$. 
\end{proof}

Let $\tilde{\sigma}=\min \{ {\sigma}^{n/2},{\sigma}_1 \}$. Until now we proved that, if the 
orbit of $P$ visits the set $W$, then when
it leaves the set 
$\tilde{W}$, it has accumulated exponential growth to the
unstable vectors 
by a
factor $\sigma^{1/2}\geq\tilde{\sigma}$ as we have seen in 
lemma
\ref{expansao}. 
As shown above,  iterates outside
$\tilde{W}$ 
contribute as well 
with the same factor.
Similar calculations hold 
for the stable vectors, and we define $\tilde{\lambda}$ through the analogous estimates for $\Phi^{-1}$. We remark that in particular these estimates improve significantly on 
the growth estimates of \cite{Rio01} which already imply the existence of a 
non-uniformly hyperbolic splitting as defined in the introduction.

\begin{lemma} 
\label{sequencia}
If $z$ does not belong to the orbit of homoclinic tangency, 
$v\in E^u_z$, and $w\in E^s_z$, then there exist sequences
 $n_k\rightarrow \infty$ and $n_l\rightarrow \infty$, and a 
positive number $c(z)$ such that 
$\Vert D\F^{n_k}_zv\Vert > c(z)\tilde{\sigma}^{n_k}\Vert v \Vert$ 
and $\Vert D\F^{-n_l}_zw\Vert > c(z)\tilde{\lambda}^{-n_k}
\Vert w \Vert$. 

In particular, $limsup_{n\rightarrow \infty}
\frac{\log{\Vert D\F^n_z\Vert}}{n}\geq \log{\tilde{\sigma}}$ 
and $liminf_{n\rightarrow \infty}\frac{\log{\Vert D\F^n_z\Vert}}{n}
\leq \log{\tilde{\lambda}}.$ 
\end{lemma} 

\begin{proof} Let $z$ be a point in $\Lambda$ outside the orbit of homoclinic tangency. For $i\in\mathbb{N}$, 
let $z_i=\Phi^i(z)$ and $v_i=D\F^i_zv$. Let $N(z)=\{ i\in 
\mathbb{N}:z_i\in W \}$. If $N(z)$ is 
non-empty and infinite, let $N(z)=\{ i_1, i_2, i_3, \ldots \}$ and 
consider $n_k$ as the smaller natural number bigger than $i_k$ such 
that $z_{n_k}\notin R_1$. In this case, we have $i_k<n_k<i_{k+1}$, 
for any $k\in\mathbb{Z}$. As a consequence of the estimates in lemmas (\ref{expansao}) and (\ref{first}) , we have that 

$$\Vert v_{n_k} \Vert > \tilde{\sigma}^{n_k-i_k}\Vert v_{i_k} 
\Vert,$$  
and hence 
$$\Vert v_{n_k} \Vert > \tilde{\sigma}^{n_k}\Vert v_{i_1} \Vert.$$ 
Taking $c(z)=\Vert v \Vert /\Vert v_{i_1}\Vert$, the result follows 
in this case. 

If $N(z)$ is finite, $i_l=\sup N(z)$, then for all $n\in \mathbb{N}$ larger 
than $i_l$, we have $\Vert v_i \Vert > \tilde{\sigma}^{i-i_l} \Vert 
v_{i_l}\Vert$, and the lemma is true for $z$, choosing $c(z)=\Vert v 
\Vert /\Vert v_{i_l}\Vert$. 

If $N(z)$ is empty, then, there are three cases to consider. If the orbit 
of $z$ is not in $\tilde{W}$, then we have exponential growth beginning at the first iterate. If $z\in \tilde{W}$ and leaves it eventually, we wait until it does so to have exponential growth, and argue as before. If $z_i\in \tilde{W}$ and never leaves it, then $z_i\in R_1$ for all 
$i\in\mathbb{N}$, meaning that $z$ is in the local stable manifold of 
$P=(0,0)$. In this case, any non-horizontal vector based in $z$ will eventually 
grow exponentially by a factor $\sigma-\varepsilon >\tilde{\sigma}$. 

Similar arguments prove 
the statement for the stable vectors.
\end{proof}

Now, define $\Lambda'=\Lambda \setminus \{ \Phi^n(q,o): n\in \mathbb Z \}$.
\begin{lemma}\label{medida}
If $\mu$ is an invariant probability measure supported 
on $\Lambda$, then $\mu(\Lambda')=1$
\end{lemma}

\begin{proof}
Let $q$ be a point in the orbit of homoclinic tangency, and suppose 
that $\mu(\{ q_i,\;i\in\mathbb{Z} \} )>0$. Then, for some $i_0$, we 
must have $\mu\{ q_{i_0} \} >0$, and the entire orbit would have 
infinite measure. 
\end{proof}

\begin{proof}[Proof of the theorem] 
Let $\mu$ be an invariant measure, and
$B$ the subset of $\Lambda'$ for which the lyapunov exponents
exist. Then, by the Oseledet's theorem 
$\mu( \mathcal B \cap \Lambda')=1$. By lemma (\ref{sequencia})
 all Lyapunov
exponents are outside the interval $(\log{\tilde{\lambda}}, \log{\tilde{\sigma}})$, 
for all $z\in B$.  
\end{proof}

\end{document}